\documentclass{icm2010}

\def\sqr#1#2{{\vcenter{\vbox{\hrule height.#2pt
              \hbox{\vrule width.#2pt height#1pt \kern#1pt \vrule width.#2pt}
              \hrule height.#2pt}}}}
\def\signed #1{{\unskip\nobreak\hfil\penalty50
              \hskip2em\hbox{}\nobreak\hfil#1
              \parfillskip=0pt \finalhyphendemerits=0 \par}}
\def\endpf{\signed {$\sqr69$}}

\def\dbR{{\mathop{\rm l\negthinspace R}}}

\def\3n{\negthinspace \negthinspace \negthinspace }
\def\2n{\negthinspace \negthinspace }
\def\1n{\negthinspace }

\def\ds{\displaystyle}

\def\dbR{{\mathop{\rm l\negthinspace R}}}
\def\={\buildrel \triangle \over =}

\def\resp{{\it resp. }}
\def\a{\alpha}
\def\b{\beta}

\def\e{\varepsilon}

\def\l{\lambda}

 \def\n{\nabla}

\def\t{\times}
\def\f{\varphi}
\def\th{\theta}
\def\o{\omega}

\def\i{\infty}
\def\ns{\noalign{\ss} }

\def\pa{\partial}
\def\G{\Gamma}
\def\D{\Delta}

\def\Si{\Sigma}

\def\O{\Omega}

\def\cA{{\cal A}}

\def\cC{{\cal C}}

\def\cF{{\cal F}}

\def\cO{{\cal O}}
\def\cP{{\cal P}}

\def\cl{{\cal l}}
\def\ov{{\overline v}}

\def\ss{\smallskip}
\def\ms{\medskip}

\def\q{\quad}
\def\qq{\qquad}
\def\hb{\hbox}

\def\lan{\mathop{\langle}}
\def\ran{\mathop{\rangle}}

\def\max{\mathop{\rm max}}
\def\min{\mathop{\rm min}}
\def\exp{\mathop{\rm exp}}
\def\sup{\mathop{\rm sup}}

\def\pa{\partial}

\def\cd{\cdot}

\def\ae{\hbox{\rm a.e.{ }}}
\def\as{\hbox{\rm a.s.{ }}}

\def\supp{\hbox{\rm supp$\,$}}

\def\cl{\overline}

\def\|{\Big |}
\def\({\Big (}
\def\){\Big )}
\def\[{\Big[}
\def\]{\Big]}
\def\be{\begin{equation}}
\def\bel{\begin{equation}\label}
\def\ee{\end{equation}}
\def\bt{\begin{theorem}}
\def\bcd{\begin{condition}}
\def\ecd{\end{condition}}
\def\et{\end{theorem}}
\def\bc{\begin{corollary}}
\def\ec{\end{corollary}}
\def\bde{\begin{definition}}
\def\ede{\end{definition}}
\def\bl{\begin{lemma}}
\def\el{\end{lemma}}
\def\bp{\begin{proposition}}
\def\ep{\end{proposition}}
\def\br{\begin{remark}}
\def\er{\end{remark}}
\def\ba{\begin{array}}
\def\ea{\end{array}}
\def\ed{\end{document}}
\def\ns{\noalign{\ms}}
\def\ds{\displaystyle}

\def\square#1{\vbox{\hrule\hbox{\vrule height#1%
     \kern#1\vrule}\hrule}}
\def\rectangle#1#2{\vbox{\hrule\hbox{\vrule height#1%
     \kern#2\vrule}\hrule}}
\def\qed{\hfill \vrule height7pt width3pt depth0pt}

\title[A unified infinite-dimensional controllability/observability theory]{A unified controllability/observability theory for some stochastic and deterministic partial differential
equations}
\author[X. Zhang]{Xu Zhang
\thanks{This work is supported by the
NSFC under grants 10831007, 60821091 and 60974035,  and the project
MTM2008-03541 of the Spanish Ministry of Science and Innovation.}%
} \contact[xuzhang@amss.ac.cn]{School of Mathematics, Sichuan
University, Chengdu 610064, China; and Key Laboratory of Systems
Control, Academy of Mathematics and Systems Science, Chinese Academy
of Sciences, Beijing 100190, China.}

\newtheorem{theorem}{Theorem}[section]
\newtheorem{prop}[theorem]{Proposition}
\newtheorem{lemma}[theorem]{Lemma}
\newtheorem{remark}[theorem]{Remark}
\newtheorem{assum}[theorem]{Condition}

\theoremstyle{definition}
\newtheorem{definition}{Definition}
\numberwithin{equation}{section}

\begin{document}

\begin{abstract} The purpose of this paper is to present a universal approach to the study of
controllability/observability problems for infinite dimensional
systems governed by some stochastic/deterministic partial
differential equations. The crucial analytic tool is a class of
fundamental weighted identities for stochastic/deterministic partial
differential operators, via which one can derive the desired global
Carleman estimates. This method can also give a unified treatment of
the stabilization, global unique continuation, and inverse problems
for some stochastic/deterministic partial differential equations.
\end{abstract}

\begin{classification}
Primary 93B05; Secondary 35Q93, 93B07.
\end{classification}

\begin{keywords} Controllability,
observability, parabolic equations, hyperbolic equations, weighted
identity.

\end{keywords}

\maketitle


\section{Introduction}

We begin with the following controlled system governed by a linear
Ordinary Differential Equation (ODE for short):
 \bel{ols1}
 \left\{
 \ba{ll}
 \ds\frac{dy(t)}{dt} =Ay(t)+Bu(t),\qq t>0,\\
  \ns
 y(0)=y_0.
 \ea\right.
 \ee
In (\ref{ols1}), $A\in \mathbb{R}^{n\times n}$, $B\in
\mathbb{R}^{n\times m}$ ($n, m\in\mathbb{N}$), $y(\cd)$ is the {\it
state variable}, $u(\cd)$ is the {\it control variable},
$\mathbb{R}^n$ and $\mathbb{R}^m$ are the {\it state space and
control space}, respectively. System (\ref{ols1}) is said to be
exactly controllable at a time $T>0$ if for any initial state
$y_0\in\mathbb{R}^n$ and any final state $y_1\in\mathbb{R}^n$, there
is a control $u(\cd)\in L^2(0,T;\mathbb{R}^m)$ such that the
solution $y(\cd)$ of (\ref{ols1}) satisfies $ y(T)=y_1$.

The above definition of controllability can be easily extended to
abstract evolution equations. In the general setting, it may happen
that the requirement $y(T)=y_1$ has to be relaxed in one way or
another. This leads to the approximate controllability, null
controllability, and partial controllability, etc. Roughly speaking,
the controllability problem for an evolution process is driving the
state of the system to a prescribed final target state (exactly or
in some approximate way) at a finite time. Also, the above $B$ can
be unbounded for general controlled systems.

The controllability/observability theory for finite dimensional
linear systems was introduced by R.E. Kalman (\cite{15}). It is by
now the basis of the whole control theory. Note that a finite
dimensional system is usually an approximation of some infinite
dimensional system. Therefore, stimulated by Kalman's work, many
mathematicians devoted to extend it to more general systems
including infinite dimensional systems, and its nonlinear and
stochastic counterparts. However, compared with Kalman's classical
theory, the extended theories are not very mature.

Let us review rapidly the main results of Kalman's theory. First of
all, it is shown that: {\it System (\ref{ols1}) is exactly
controllable at a time $T$ if and only if} $ \hbox{\rm
rank$\,$}[B,AB, \cdots, A^{n-1} B] = n$. However, this criterion is
not applicable for general infinite dimensional systems. Instead, in
the general setting, one uses another method which reduces the
controllability problem for a controlled system to an observability
problem for its dual system. The dual system of (\ref{ols1}) reads:
 \bel{c1s1e3---1}
 \left\{
 \ba{ll}
 \ds\frac{dw}{dt}=-A^*w,\qq t\in (0,T),\\
 \ns
 w(T)=z_0.
 \ea\right.
 \ee
It is shown that: {\it System (\ref{ols1}) is exactly controllable
at some time $T$ if and only if the following observability
inequality (or estimate) holds}
 \bel{c1s1e6}
 |z_0|^2\le C\int_0^T\left|B^*w(t)\right|^2dt,\q
 \forall\;z_0\in\mathbb{R}^n.
 \ee
Here and henceforth, $C$ denotes a generic positive constant, which
may be different from one place to another. We remark that similar
results remain true in the infinite dimensional setting, where the
theme of the controllability/observability theory is to establish
suitable observability estimates through various approaches.

Systems governed by Partial Differential Equations (PDEs for short)
are typically infinite dimensional. There exists many works on
controllability/observability of PDEs. Contributions by D.L.~Russell
(\cite{Russell}) and by J.L.~Lions (\cite{Lions}) are classical in
this field. In particular, since it stimulated many in-depth
researches on related problems in PDEs, J.L. Lions's paper
\cite{Lions} triggered extensive works addressing the
controllability/observability of infinite dimensional controlled
system. After \cite{Lions}, important works in this field can be
found in
\cite{AI,BLR,Coron,FGIP,Fu1,FI,KoL,LR,TTLi,TW,Yama,Zua1,Zua2}. For
other related works, we refer to \cite{Isakov,LY} and so on.

The controllability/observability of PDEs depends strongly on the
nature of the underlying system, such as time reversibility or not,
and propagation speed of solutions, etc. The wave equation and the
heat equation are typical examples. Now it is clear that essential
differences exist between the controllability/observability theories
for these two equations. Naturally, one expects to know whether some
relationship exist between the controllability/observability
theories for these two equations of different nature. Especially, it
would be quite interesting to establish, in some sense and to some
extend, a unified controllability/observability theory for parabolic
equations and hyperbolic equations. This problem was initially
studied by D.L. Russell (\cite{28}).

The main purpose of this paper is to present the author's and his
collaborators' works with an effort towards a unified
controllability/observability theory for stochastic/deterministic
PDEs. The crucial analytic tool we employ is a class of elementary
pointwise weighted identities for partial differential operators.
Starting from these identities, we develop a unified approach, based
on global Carleman estimate. This universal approach not only
deduces the known controllability/observability results (that have
been derived before via Carleman estimates) for the linear
parabolic, hyperbolic, Schr\"odinger and plate equations, but also
provides new/sharp results on controllability/observability, global
unique continuation, stabilization and inverse problems for some
stochastic/deterministic linear/nonlinear PDEs.

The rest of this paper is organized as follows. Section \ref{s2}
analyzes the main differences between the existing
controllability/observability theories for parabolic equations and
hyperbolic equations. Sections \ref{s3} and \ref{s4} address, among
others, the unified treatment of the controllability/observability
problem for deterministic PDEs and stochastic PDEs, respectively.

\section{Main differences between the known theories}\label{s2}
In the sequel, unless otherwise indicated, $G$ stands for a bounded
domain (in $\mathbb{R}^n$) with a boundary $\G\in C^2$, $G_0$
denotes an open non-empty subset of $G$, and $T$ is a given positive
number. Put $Q=(0,T)\times G$, $Q_{G_0}=(0,T)\times G_0$ and
$\Si=(0,T)\times\G$.

We begin with a controlled heat equation:
 \be\label{he}\left\{\ba{ll} y_{t}-\Delta y=\chi_{G_0}(x)u(t,x)
 \quad & \hbox{in }Q,\\
 y=0 \quad & \hbox{on }\Si,\\
 y(0)=y_{0} \quad & \hbox{in }G \ea\right.
 \ee
and a controlled wave equation:
 \be\label{we}\left\{\ba{ll} y_{tt}-\D y=\chi_{G_0}(x)u(t,x) \quad &
 \hbox{in }Q,\\
 y=0 \quad & \hbox{on }\Si,\\
 y(0)=y_{0}, \q y_{t}(0)=y_{1} \quad & \hbox{in }G.
 \ea\right.\ee
In (\ref{he}), $y$ and $u$ are the {\it state variable} and {\it
control variable}, the {\it state space} and {\it control space} are
chosen to be $L^{2}(G)$ and $L^{2}(Q_{G_0})$, respectively; while in
(\ref{we}), $(y,y_t)$ and $u$ are the {\it state variable} and {\it
control variable}, $H^{1}_{0}(G)\times L^{2}(G)$ and
$L^{2}(Q_{G_0})$ are respectively the {\it state space} and {\it
control space}. System (\ref{he}) is said to be null controllable
(\resp approximately controllable) in $L^2(G)$ if for any given
$y_0\in L^2(G)$ (\resp for any given $\e>0$, $y_0, y_1\in L^2(G)$),
one can find a control $u\in L^2(Q_{G_0})$ such that the weak
solution $y(\cd)\in C([0,T];L^2(G))\cap C((0,T];H_0^1(G))$ of
(\ref{he}) satisfies $y(T)=0$ (\resp $|y(T)-y_1|_{L^2(G)}\le \e$).
In the case of null controllability, the corresponding control $u$
is called a null-control (with initial state $y_0$). Note that, due
to the smoothing effect of solutions to the heat equation, the exact
controllability for (\ref{he}) is impossible, i.e., the above $\e$
cannot be zero. On the other hand, since one can rewrite system
(\ref{we}) as an evolution equation in a form like (\ref{ols1}), it
is easy to define the exact controllability of this system. The dual
systems of (\ref{he}) and (\ref{we}) read respectively
 \begin{eqnarray}\label{olwsystem1}
\left\{
\begin{array}{lll}
\ds  \psi_t+ \D \psi  = 0
&\mbox{ in } Q,\\
\psi = 0 & \mbox{ on } \Si, \\
\psi(T) = \psi_0&\mbox{ in } G
\end{array}
\right.
\end{eqnarray}
 and
\begin{eqnarray}\label{olhsystem1}
\left\{
\begin{array}{lll}
\ds  \psi_{tt} - \D \psi  = 0
&\mbox{ in } Q,\\
\psi = 0 & \mbox{ on } \Si, \\
\psi(T) = \psi_0, \q\psi_t(T) = \psi_1 &\mbox{ in } G.
\end{array}
\right.
\end{eqnarray}

The controllability/observability theories for parabolic equations
and hyperbolic equations turns out to be quite different. First of
all, we recall the related result for the heat equation.

\bt\label{2t1}
{\rm (\cite{LR})} Let $G$ be a bounded domain of class $C^{\infty}$.
Then: i) System (\ref{he}) is null controllable and approximately
controllable in $L^2(G)$ at time $T$; ii) Solutions of equation
(\ref{olwsystem1}) satisfy
 \bel{icm1}
|\psi(0)|_{L^2(G)}\le C|\psi|_{L^2(Q_{G_0})},\qq\forall\;\psi_0\in
L^2(G).
 \ee
\et

Since solutions to the heat equation have an infinite propagation
speed, the ``waiting" time $T$ can be chosen as small as one likes,
and the control domain $G_0$ dose not need to satisfy any geometric
condition but being open and non-empty. On the other hand, due to
the time irreversibility and the strong dissipativity of
(\ref{olwsystem1}), one cannot replace $|\psi(0)|_{L^2(G)}$ in
inequality (\ref{icm1}) by $|\psi_0|_{L^2(G)}$.

Denote by $\{\mu_i\}^\infty_{i=1}$ the eigenvalues of the homogenous
Dirichlet Laplacian on $G$, and $\{\f_i\}^{\infty}_{i=1}$ the
corresponding eigenvectors satisfying $|\f_i|_{L^2(G)} = 1$. The
proof of Theorem \ref{2t1} is based on the following observability
estimate on sums of eigenfunctions for the Laplacian (\cite{LR}):

\bt\label{icm2t1}
Under the assumption of Theorem \ref{2t1}, for any $r
> 0$, it holds
\begin{equation}\label{lr}
\sum_{\mu_{i} \leq r}|a_{i}|^{2} \leq C
e^{C\sqrt{r}}{\int_{G_0}}\bigg| \sum_{\mu_{i} \leq r}a_{i}\f_{i}(x)
\bigg|^2dx,\ \ \forall\, \{a_{i}\}_{\mu_{i} \leq r }\hb{ with }a_{i}
\in \mathbb{C}.
\end{equation}
\et

Note that Theorem \ref{icm2t1} has some other applications in
control problems of PDEs (\cite{Lopez,Miller,Wang,33,Zua1,Zua2}).
Besides, to prove Theorem \ref{2t1}, one needs to utilize a time
iteration method (\cite{LR}), which uses essentially the Fourier
decomposition of solutions to (\ref{olwsystem1}) and especially, the
strong dissipativity of this equation. Hence, this method cannot be
applied to conservative systems (say, system (\ref{we})) or the
system that the underlined equation is time-dependent.

As for the controllability/observability for the wave equation, we
need to introduce the following notations. Fix any $x_0\in
\mathbb{R}^n$, put
 \be\label{boundary}
 \G_0\=\big\{x\in \G\;\big|\; (x-x_0)\cdot\nu(x)>0\big\},
 \ee
where $\nu(x)$ is the unit outward normal vector of $G$ at
$x\in\Gamma$. For any set $S\in \mathbb{R}^n$ and $\e>0$, put ${\cal
O}_{\epsilon}(S)=\big\{y\in \mathbb{R}^n\ \big|\
 |y-x|<\e\hbox { for some }x\in S\big\}$.

The exact controllability of system (\ref{we}) is equivalent to the
following {\it observability estimate} for system
(\ref{olhsystem1}):
\begin{equation}\label{oinlhsystem1}
 |(\psi_0,\psi_1)|_{L^2(G)\t H^{-1}(G)} \1n\leq\1n
C|\psi|_{L^2(Q_{G_0})}, \ \ \ \forall\, (\psi_0,\psi_1)\in L^2(G)\t
H^{-1}(G).
\end{equation}
Note that the left hand side of (\ref{oinlhsystem1}) can be replaced
by $|(\psi(0),\psi_t(0))|^2_{L^2(G)\t H^{-1}(G)}$ (because
(\ref{olhsystem1}) is conservative). The following classical result
can be found in \cite{Lions}.

\bt\label{thlr} Assume $G_0=O_{\e}(\G_0)\cap G$ and $\ds
T_0=2\sup_{x\in G\setminus G_0}|x-x_0|$. Then, inequality
(\ref{oinlhsystem1}) holds for any time $T>T_0$.
\et

The proof of Theorem \ref{thlr} is based on a classical Rellich-type
multiplier method. Indeed, it is a consequence of the following
identity (e.g. \cite{31}):

\begin{prop}
Let $h\=(h^1, \cdots, h^n):\
\mathbb{R}\t\mathbb{R}^{n}\to\mathbb{R}^n$ be a vector field of
class $C^1$. Then for any $z\in C^2(\mathbb{R}\t\mathbb{R}^n)$, it
holds that
$$
\ba{ll}
\ds\n\cd\Big\{2(h\cd\n z)(\n z)+h\Big[z_t^2-\sum_{i=1}^nz_{x_i}^2\Big]\Big\}\\
 \ns\ds
=-2(z_{tt}-\D z)h\cd\n z+(2z_th\cd\n z)_t-2z_th_t\cd\n z\\
 \ns\ds\q+(\n\cd
h)\Big[z_t^2-\sum_{i=1}^nz_{x_i}^2\Big]
+2\sum_{i,j=1}^n\Big(\frac{\pa h^j}{\pa x_i}z_{x_i}z_{x_j}\Big).
 \ea
$$
\end{prop}

The observability time $T$ in Theorem \ref{thlr} should be large
enough. This is due to the finite propagation speed of solutions to
the wave equation (except when the control is acting in the whole
domain $G$). On the other hand, it is shown in \cite{BLR} that exact
controllability of (\ref{we}) is impossible without geometric
conditions on $G_0$. Note also that, the multiplier method rarely
provides the optimal control/observation domain and minimal
controll/observation time except for some very special geometries.
These restrictions are weakened by the microlocal analysis
(\cite{BLR}). In \cite{BLR,Burq,BG}, the authors proved that,
roughly speaking, inequality (\ref{oinlhsystem1}) holds if and only
if every ray of Geometric Optics that propagates in $G$ and is
reflected on its boundary $\G$ enters $G_0$ at time less than $T$.

The above discussion indicates that the results and methods for the
controllability/observability of the heat equation differ from those
of the wave equation. As we mentioned before, this leads to the
problem of establishing a unified theory for the
controllability/observability of parabolic equations and hyperbolic
equations. The first result in this direction was given in
\cite{28}, which showed that {\it the exact controllability of the
wave equation implies the null controllability of the heat equation
with the same controller but in a short time}. Further results were
obtained in \cite{Lopez, 33}, in which organic connections were
established for the controllability theories between parabolic
equations and hyperbolic equations. More precisely, it has been
shown that: {\it i) By taking the singular limit of some exactly
controllable hyperbolic equations, one gives the null
controllability of some parabolic equations (\cite{Lopez}); and ii)
Controllability results of the heat equation can be derived from the
exact controllability of some hyperbolic equations (\cite{33})}.
Other interesting related works can be found in
\cite{Miller,Phung,TW}. In the sequel, we shall focus mainly on a
unified treatment of the controllability/observability for both
deterministic PDEs and stochastic PDEs, from the methodology point
of view.

\section{The deterministic case}\label{s3}

The key to solve controllability/observability problems  for PDEs is
the obtention of suitable observability inequalities for the
underlying homogeneous systems. Nevertheless, as we see in Section
\ref{s2}, the techniques that have been developed to obtain such
estimates depend heavily on the nature of the equations, especially
when one expects to obtain sharp results for time-invariant
equations. As for the time-variant case, in principle one needs to
employ Carleman estimates, see \cite{FI} for the parabolic equation
and \cite{31} for the hyperbolic equation. The Carleman estimate is
simply a weighted energy method. However, at least formally, the
Carleman estimate used to derive the observability inequality for
parabolic equations is quite different from that for hyperbolic
ones. The main purpose of this section is to present a universal
approach for the controllability/observability of some deterministic
PDEs. Our approach is based on global Carleman estimates via a
fundamental pointwise weighted identity for partial differential
operators of second order (It was established in \cite{Fu1,fx}. See
\cite{[20]} for an earlier result). This approach is stimulated by
\cite{LRS,KK}, both of which are addressed for ill-posed problems.

\subsection{A stimulating example}

The basic idea of Carleman estimates is available in proving the
stability of ODEs (\cite{[20]}). Indeed, consider an ODE in
$\mathbb{R}^{n}$:
 \be\label{ODE}
 \left\{\ba{ll} x_t(t)=a(t)x(t),\quad
 t\in [0,T],
 \\x(0)=x_{0}, \ea
 \right.
 \ee
where $a\in L^{\i}(0,T)$. A well-known simple result reads: {\it
Solutions of (\ref{ODE}) satisfy
 \bel{3e1}
 \max_{t\in[0,T]}|x(t)|\le C|x_0|,\qq\forall\;x_0\in \mathbb{R}^n.
 \ee}

{\it A Carleman-type Proof of (\ref{3e1})}. For any $\l\in
\mathbb{R}$, by (\ref{ODE}), one obtains
 \be
 \label{CD}
 \frac{d}{dt}\(e^{-\l t}|x(t)|^{2}\)=-\l e^{-\l t}|x(t)|^{2}+2e^{-\l t} x_t(t)\cd x(t)
 =(2a(t)-\l)e^{-\l t}|x(t)|^{2}.
 \ee
Choosing $\l$ large enough so that $2a(t)-\l\le 0$ for a.e. $t\in
(0,T)$, we find that
 $$
 |x(t)|\leq e^{\l T/2}|x_{0}|,\quad t\in[0,T],
 $$
which proves (\ref{3e1}). \endpf

\br
By  (\ref{CD}), we see the following pointwise identity:
 \bel{3e2}
 2e^{-\l t} x_t(t)\cd x(t)=\frac{d}{dt}\(e^{-\l t}|x(t)|^{2}\)+\l e^{-\l t}|x(t)|^{2}.
 \ee
Note that  $ x_t(t)$ is the principal operator of the first equation
in (\ref{ODE}). The main idea of (\ref{3e2}) is to establish a
pointwise identity (and/or estimate) on the principal operator $
x_t(t)$ in terms of the sum of a ``divergence" term
$\frac{d}{dt}(e^{-\l t}|x(t)|^{2})$ and an ``energy" term $\l e^{-\l
t}|x(t)|^{2}$. As we see in the above proof, one chooses $\l$ to be
big enough to absorb the undesired terms. This is the key of all
Carleman-type estimates. In the sequel, we use exactly the same
method, i.e., the method of Carleman estimate via pointwise
estimate, to derive observability inequalities for both parabolic
equations and hyperbolic equations.
\er

\subsection{Pointwise weighted identity}

We now show a fundamental pointwise weighted identity for general
partial differential operator of second order. Fix real functions
$\a,\;\b\in C^1(\mathbb{R}^{1+m})$ and $b^{jk}\in
C^1(\mathbb{R}^{1+m})$ satisfying $
 b^{jk}=b^{kj}$ ($j,k=1,2,\cdots,m$).
Define a formal differential operator of second order: $\ds \cP
z\=(\a+i\b)z_t+\sum_{j,k=1}^m\big(b^{jk}z_{x_j}\big)_{x_k}$,
$i=\sqrt{-1}$. The following identity was established in
\cite{Fu1,fx}:

\bt\label{uTheorem}
Let $z\in C^2(\mathbb{R}^{1+m}; \;\mathbb{C})$ and $\ell\in
C^2(\mathbb{R}^{1+m};\mathbb{R})$. Put $\th=e^\ell$ and $v=\th z$.
Let $a, b,\l\in\mathbb{R}$ be parameters. Then
 \bel{c2a1}\ba{ll}\ds
 \th(\cP z\overline {I_1}+\overline{\cP z} I_1)+M_t+\sum_{k=1}^m\pa_{x_k}V^k\\
\ns\ds=2|I_1|^2+\sum_{j,k,j',k'=1}^m\[2(b^{j'k}\ell_{x_{j'}})_{x_{k'}}b^{jk'}-(b^{jk}b^{j'k'}\ell_{x_{j'}})_{x_{k'}}+\frac{1}{2}(\a
b^{jk})_t\\
\ns\ds\q
-ab^{jk}b^{j'k'}\ell_{x_{j'}x_{k'}}\](v_{x_k}\ov_{x_j}+\ov_{x_k}
v_{x_j})+\[-
\sum_{j,k=1}^mb^{jk}_{x_k}\ell_{x_j}+b\l\](I_1\ov+\overline
 {I_1}v)\\
\ns\ds\q+i\sum_{j,k=1}^m\Big\{[(\b
b^{jk}\ell_{x_j})_t+b^{jk}(\b\ell_t)_{x_j}](\ov_{x_k}v-v_{x_k}\ov)\\
\ns\ds\q+[(\b b^{jk}\ell_{x_j})_{x_k}+a\b b^{jk}\ell_{x_jx_k}](\ov
v_t-v\ov_t)\Big\}-\sum_{j,k=1}^mb^{jk}\a_{x_k}(v_{x_j}\ov_t+\ov_{x_j}v_t)\\
\ns\ds\q-a\sum_{j,k,j',k'=1}^mb^{jk}(b^{j'k'}\ell_{x_{j'}x_{k'}})_{x_k}(\ov_{x_j}v+v_{x_j}\ov)+B|v|^2,
 \ea\ee
where
 $$
 \left\{
 \ba{ll}
\ds I_1\=i\b v_t-\a\ell_tv+\sum_{j,k=1}^m(b^{jk}v_{x_j})_{x_k}+Av,\\
\ns \ds
A\=\sum_{j,k=1}^mb^{jk}\ell_{x_j}\ell_{x_k}-(1+a)\sum_{j,k=1}^mb^{jk}\ell_{x_jx_k}-b\l,
 \\
\ns \ds B\=(\a^2\ell_t+\b^2\ell_{t}-\a A)_t\\ \ns\ds
\qq+2\sum_{j,k=1}^m\[(b^{jk}\ell_{x_j}A)_{x_k}-(\a
b^{jk}\ell_{x_j}\ell_{t})_{x_k}+a(A-\a\ell_t)b^{jk}\ell_{{x_j}{x_k}}\],
\\
\ns \ds M\=\[(\a^2+\b^2)\ell_t-\a A\]
|v|^2+\a\sum_{j,k=1}^mb^{jk}v_{x_j}\ov_{x_k}\\ \ns\ds
\qq+i\b\sum_{j,k=1}^mb^{jk}\ell_{x_j}(\ov_{x_k}v-v_{x_k}\ov),
\\
\ns \ds  V^k\=\sum_{j,j',k'=1}^m\Big\{-i
\b\[b^{jk}\ell_{x_j}(v\ov_t-\ov
 v_t)+b^{jk}\ell_t(v_{x_j}\ov-\ov_{x_j}v)\]\\
\ns \ds\qq\qq\qq\q -\a b^{jk}(v_{x_j}\ov_t+\ov_{x_j}v_t)\\
\ns \ds\qq\qq\qq\q+(2b^{jk'}b^{j'k}-b^{jk}b^{j'k'})\ell_{x_j}(v_{x_{j'}}\ov_{x_{k'}}+\ov_{x_{j'}}v_{x_{k'}}) \\
\ns \ds\qq\qq\qq\q -ab^{j'k'}\ell_{x_{j'}x_{k'}}
b^{jk}(v_{x_j}\ov+\ov_{x_j}v)+2b^{jk}(A\ell_{x_j}-\a\ell_{x_j}\ell_t)|v|^2\Big\}.
 \ea\right.
$$ \et

As we shall see later, Theorem \ref{uTheorem} can be applied to
study the controllability/observability as well as the stabilization
of parabolic equations and hyperbolic equations. Also, as pointed by
\cite{Fu1}, starting from Theorem \ref{uTheorem}, one can deduce the
controllability/observability for the Schr\"odinger equation and
plate equation appeared in \cite{LTZ2} and \cite{32}, respectively.
Note also that, Theorem \ref{uTheorem} can be applied to study the
controllability of the linear/nonlinear complex Ginzburg-Landau
equation (see \cite{Fu1,fx,[26]}).

\subsection{Controllability/Observability of Linear PDEs}
In this subsection, we show that, starting from Theorem
\ref{uTheorem}, one can establish sharp
observability/controllability results for both parabolic systems and
hyperbolic systems.

We need to introduce the following assumptions.

\begin{assum}\label{condition of d-0} Matrix-valued function
$\big(p^{ij}\big)_{1\le i,j\le n}\in C^{1}(\overline{Q};
\mathbb{R}^{n\times n})$ is uniformly positive definite.
\end{assum}

\begin{assum}\label{condition of d-1} Matrix-valued function
$\big(h^{ij}\big)_{1\le i,j\le n}\in C^{1}(\cl{G};
\mathbb{R}^{n\times n})$ is uniformly positive definite.
\end{assum}

Also, for any $N\in \mathbb{N}$, we introduce the following

\begin{assum}
\label{condition of d-2} Matrix-valued functions  $a\in
L^{\infty}(0,T;L^{p}(G;\mathbb{R}^{N\t N}))$ for some $p \in
[n,\infty]$, and $a_1^1, \cdots, a_1^n, a_2  \in
L^{\infty}(Q;\mathbb{R}^{N\t N})$.
\end{assum}

Let us consider first the following parabolic system:
\begin{eqnarray}\label{system4}
\left\{
\begin{array}{lll}
\ds\varphi_t - \sum_{i,j=1}^n (p^{ij}\varphi_{x_i})_{x_j} =
a\varphi + \sum_{k=1}^n a_1^k \varphi_{x_k}, &\mbox{ in } Q,\\
\ns\ds \varphi = 0, &\mbox{ on } \Si,\\
\ns\ds \varphi(0) = \varphi^0, &\mbox{ in } G,
\end{array}
\right.
\end{eqnarray}
where $\varphi$ takes values in $\mathbb{R}^{N}$. By choosing $\a =
1$ and $\b = 0$ in Theorem \ref{uTheorem}, one obtains a weighted
identity for the parabolic operator. Along with \cite{[20]}, this
identity leads to the existing controllability/observability result
for parabolic equations (\cite{Doubova1,FI}). One can go a little
further to show the following result (\cite{DZZ}):

\begin{theorem}\label{oTheorem4}
Let Conditions \ref{condition of d-0} and \ref{condition of d-2}
hold. Then, solutions  of (\ref{system4}) satisfy
 \be\label{oest3} \2n\ba{ll}
|\varphi(T)|_{(L^2(G))^N} \\
\ns \ds\leq \exp\Big\{ C\Big[ 1+\frac{1}{T} +
T|a|_{L^{\infty}(0,T;L^{p}(G;\mathbb{R}^{N\t
N}))}+|a|_{L^{\infty}(0,T;L^{p}(G;\mathbb{R}^{N\t
N}))}^{\frac{1}{\frac{3}{2}-\frac{n}{p}}}\\
 \ns\ds \q + (1+T)\Big(\sum_{k=1}^N
|a_i^k|_{L^{\infty}(Q;\mathbb{R}^{N\t N})}\Big)^2\Big]
 \Big\}|\varphi|_{(L^2(Q_{G_0}))^N},\q \forall\; \varphi^0 \in (L^2(G))^N.
\ea \ee
\end{theorem}

Note that (\ref{oest3}) provides the observability inequality for
the parabolic system (\ref{system4}) with an explicit estimate on
the observability constant, depending on the observation time $T$,
the potential $a$ and $a_1^k$. Earlier result in this respect can be
found in \cite{Doubova1} and the references cited therein.
Inequality (\ref{oest3}) will play a key role in the study of the
null controllability problem for semilinear parabolic equations, as
we shall see later.

\begin{remark}\label{icmrem1}
It is shown in \cite{DZZ} that when $n\ge 2$, $N\ge 2$ and
$\big(p^{ij}\big)_{1\leq i,j\leq n} = I$, the exponent $\frac{2}{3}$
in $|a|_{L^{\infty}(0,T;L^{p}(G;\mathbb{R}^{N\t N}))}^{\frac{2}{3}}$
(for the case that $p=\infty$ in the inequality (\ref{oest3})) is
sharp. In \cite{DZZ}, it is also proved that the quadratic
dependence on $\ds\sum_{k=1}^N |a_i^k|_{L^{\infty}(Q;\mathbb{R}^{N\t
N})}$ is sharp under the same assumptions. However, it is not clear
whether the exponent $\frac{3}{2}-\frac{n}{p}$ in
$|a|_{L^{\infty}(0,T;L^{p}(G;\mathbb{R}^{N\t
N}))}^{\frac{1}{\frac{3}{2}-\frac{n}{p}}}$ is optimal when
$p<\infty$.
\end{remark}

Next, we consider  the following hyperbolic system:
\begin{eqnarray}\label{system5}
\left\{
\begin{array}{lll}
\ds v_{tt} - \sum_{i,j=1}^n (h^{ij}v_{x_i})_{x_j} =
av + \sum_{k=1}^n a_1^k v_{x_k} + a_2 v_t, &\mbox{ in } Q,\\
\ns\ds v = 0, &\mbox{ on } \Si,\\
\ns\ds v(0) = v^0, \q v_t(0) = v^1, &\mbox{ in } G,
\end{array}
\right.
\end{eqnarray}
where $v$ takes values in $\mathbb{R}^{N}$.

Compared with the parabolic case, one needs more assumptions on the
coefficient matrix $\big(h^{ij}\big)_{1\le i,j\le n}$ as follows
(\cite{DZZ,Fu-Yong-Zhang}):

\begin{assum}
\label{condition of d} There is a positive function $d(\cdot) \in
C^2(\overline{G})$ satisfying

i) For some constant $\mu_0 \ge 4$, it holds
 $$
 \ba{ll}\ds
\sum_{i,j=1}^n\Big\{ \sum_{i',j'=1}^n\Big[
2h^{ij'}(h^{i'j}d_{x_{i'}})_{x_{j'}} -
h^{ij}_{x_{j'}}h^{i'j'}d_{x_{i'}} \Big] \Big\}\xi^{i}\xi^{j} \geq
\mu_0
\sum_{i,j=1}^n h^{ij}\xi^{i}\xi^{j}, \\
 \ns
 \ds\qq\qq\qq\qq\qq\qq\qq\qq\qq\qq
\forall\; (x,\xi^{1},\cdots,\xi^{n}) \in \overline{G} \t
\mathbb{R}^n; \ea
 $$

ii) There is no critical point of $d(\cdot)$ in $\overline{G}$,
i.e.,
 $\ds\min_{x\in \overline{G}}|\nabla d(x)| > 0$;

iii)  $\q\ds\frac{1}{4}\sum_{i,j=1}^n h^{ij}(x)d_{x_i}(x)d_{x_j}(x)
\geq \max_{x\in\overline{G}} d(x), \q \forall x\in \overline{G}$.
\end{assum}

We put
\begin{equation}\label{tstar}
 T^* = 2\max_{x\in
\overline{G}}\sqrt{d(x)}, \qq
 \G^* \= \Big\{x\in\G \;\Big|\;
\sum_{i,j=1}^n h^{ij}(x)d_{x_i}(x)\nu_j(x)>0 \Big\}.
\end{equation}
By choosing $b^{jk}(t,x) \equiv h^{jk}(x)$ and $\a = \b = 0$ in
Theorem \ref{uTheorem} (and noting that only the symmetry condition
is assumed for $b^{jk}$ in this theorem), one obtains the
fundamental identity derived in \cite{Fu-Yong-Zhang} to establish
the controllability/observability of the general hyperbolic
equations. One can go a little further to show the following result
(\cite{DZZ}).

\begin{theorem}\label{oTheorem5}
Let Conditions \ref{condition of d-1}, \ref{condition of d-2} and
\ref{condition of d} hold, $T>T^*$ and $G_0=G\cap \cO_{\e}(\G^*)$
for some $\e>0$. Then one has the following conclusions:

\ss

1) For any $(v^0,v^1) \in (H_0^1(G))^N \t (L^2(G))^N$, the
corresponding weak solution $ v \in C([0,T];(H_0^1(G))^N)\bigcap $
$C^1([0,T];(L^2(G))^N) $ of system (\ref{system5}) satisfies
\bel{oest4} \ba{ll}
\ds |v^0|_{H_0^1(G))^N} + |v^1|_{(L^2(G))^N} \\
\ns \ds \leq \exp\Big[ C\Big(1 +
|a|_{L^{\infty}(0,T;L^{p}(G;\mathbb{R}^{N\t
N}))}^{\frac{1}{\frac{3}{2}-\frac{n}{p}}}\\
\ns\ds\q +\Big(\sum_{k=1}^N |a_i^k|_{L^{\infty}(Q;\mathbb{R}^{N\t
N})}+|a_2|_{L^{\infty}(Q;\mathbb{R}^{N\t N})}\Big)^2\Big) \Big]
\Big| \frac{\pa v}{\pa \nu} \Big|_{(L^2((0,T)\t\G^*))^N}. \ea \ee

\ss

2) If $a_1^k \equiv 0$ ($k = 1,\cdots,n$) and $a_2 \equiv 0$, then
for any $(v^0,v^1) \in (L^2(G))^N \t (H^{-1}(G))^N$, the weak
solution $ v \in C([0,T];(L^2(G))^N)\bigcap C^1([0,T];(H^{-1}(G))^N)
$ of system (\ref{system5}) satisfies
 \bel{oest5}
 \ba{ll}
 |v^0|_{(L^2(G))^N} + |v^1|_{H^{-1}(G))^N} \\
 \ns
 \ds\leq \exp\Big[ C\Big(1 +
|a|_{L^{\infty}(0,T;L^{p}(G;\mathbb{R}^{N\t
N}))}^{\frac{1}{\frac{3}{2}-\frac{n}{p}}}\Big) \Big]
|v|_{(L^2(Q_{G_0}))^N}.
 \ea
 \ee
\end{theorem}

As we shall see in the next subsection, inequality (\ref{oest5})
plays a crucial role in the study of the exact controllability
problem for semilinear hyperbolic equations.

\begin{remark}\label{ICMr1}
As in the parabolic case, it is shown in \cite{DZZ} that the
exponent $\ds\frac{2}{3}$ in the estimate
$|a|_{L^{\infty}(0,T;L^{p}(G;\mathbb{R}^{N\t N}))}^{\frac{2}{3}}$ in
(\ref{oest5}) (for the special case $p =\infty$) is sharp for $n\ge
2$ and $N\ge 2$. Also, the exponent $2$ in the term
$\ds\Big(\sum_{k=1}^N |a_i^k|_{L^{\infty}(Q;\mathbb{R}^{N\t
N})}+|a_2|_{L^{\infty}(Q;\mathbb{R}^{N\t N})}\Big)^2$ in
(\ref{oest4}) is sharp. However, it is unknown whether the estimate
is optimal for the case that $p<\infty$.
\end{remark}

By the standard duality argument, Theorems \ref{oTheorem4} and
\ref{oTheorem5} can be applied to deduce the controllability results
for parabolic systems and hyperbolic systems, respectively. We omit
the details.

\subsection{Controllability of Semi-linear PDEs}

The study of exact/null controllability problems for semi-linear
PDEs began in the 1960s. Early works in this respect were mainly
devoted to the local controllability problem. By the local
controllability of a system, we mean that the controllability
property holds under some smallness assumptions on the initial data
and/or the final target, or the Lipschitz constant of the
nonlinearity.

In this subsection we shall present some global controllability
results for both semilinear parabolic equations and hyperbolic
equations. These results can be deduced from Theorems
\ref{oTheorem4} and \ref{oTheorem5}, respectively.

Consider first the following controlled semi-linear parabolic
equation:
\begin{eqnarray}\label{slheat1}
\left\{
\begin{array}{lll}
\ds y_t - \sum_{i,j=1}^n (p^{ij}y_{x_i})_{x_j} + f(y,\nabla y) =\chi_{G_0} u, &\mbox{ in } Q,\\
\ns\ds y = 0, &\mbox{ on } \Si,\\
\ns\ds y(0) = y_0, &\mbox{ in }G.
\end{array}
\right.
\end{eqnarray}
For system (\ref{slheat1}), the state variable and control variable,
state  space and control space, controllability, are chosen/defined
in a similar way as for system (\ref{he}). Concerning the
nonlinearity $f(\cd,\cd)$, we introduce the following assumption
(\cite{Doubova1}).

\begin{assum}
\label{condition of d4}
Function $f(\cdot,\cdot) \in C(\mathbb{R}^{1+n})$ is locally
Lipschitz-continuous. It satisfies $f(0,0) = 0$ and
\begin{eqnarray}\label{f1}
\left\{
\begin{array}{ll}
 \ds \lim_{|(s,p)|\to\infty}\frac{\int_0^1
f_s (\tau s, \tau p)d\tau}{\ln^{\frac{3}{2}}(1+ |s| + |p|)} =
0,\\
\ns\ds \lim_{|(s,p)|\to\infty}\frac{|(\int_0^1 f_{p_1} (\tau s, \tau
p)d\tau,\cdots,\int_0^1 f_{p_n} (\tau s, \tau
p)d\tau)|}{\ln^{\frac{1}{2}}(1+ |s| + |p|)} = 0,
\end{array}
\right.
\end{eqnarray}
where $p = (p_1,\cdots, p_n)$.
\end{assum}

As shown in \cite{Doubova1} (See \cite{barbu} and the references
therein for earlier results), linearizing the equation, estimating
the cost of the control in terms of the size of the potential
entering in the system (thanks to Theorem \ref{oTheorem4}), and
using the classical fixed point argument, one can show the following
result.

\begin{theorem}\label{3.1}
Assume that Conditions \ref{condition of d-0} and \ref{condition of
d4} hold. Then system (\ref{slheat1}) is null controllable.
\end{theorem}

In particular, Theorem \ref{3.1} provides the possibility of
controlling some blowing-up equations. More precisely, assume that
$f(s,p)\equiv f(s)$ in system (\ref{slheat1}) has the form
\begin{equation}\label{eq4connections}
f(s)=-s\ln^r(1+|s|),\qq r\ge 0.
\end{equation}
When $r>1$, solutions of (\ref{slheat1}), in the absence of control,
i.e. with $u\equiv0$, blow-up in finite time. According to Theorem
\ref{3.1} the process can be controlled, and, in particular, the
blow-up can be avoided when $1<r\le 3/2$. By the contrary, it is
proved in \cite{barbu,FZ} that for some nonlinearities $f$
satisfying
 \begin{equation}\label{zzeq3connections}
\ds\lim_{\mid s\mid\to\infty}\frac{\mid f(s)\mid}{s\ln^r(1+\mid
s\mid)}=0,
\end{equation}
where $r>2$, the corresponding system is not controllable. The
reason is that the controls cannot help the system to avoid blow-up.

\begin{remark}
It is still an unsolved problem whether the controllability holds
for system (\ref{slheat1}) in which the nonlinear function $f(\cd)$
satisfies (\ref{zzeq3connections}) with $3/2\leq r\leq 2$. Note
that, the growth condition in (\ref{f1}) comes from the
observability inequality (\ref{oest3}). Indeed, the logarithmic
function in (\ref{f1}) is precisely the inverse of the exponential
one in (\ref{oest3}). According to Remark \ref{icmrem1}, the
estimate (\ref{oest3}) cannot be improved, and therefore, the usual
linearization  approach cannot lead to any improvement of the growth
condition (\ref{f1}).
\end{remark}

Next, we consider the following controlled semi-linear hyperbolic
equation:
\begin{eqnarray}\label{system3}
\left\{
\begin{array}{lll}
\ds y_{tt} - \sum_{i,j=1}^n (h^{ij}y_{x_i})_{x_j} = h(y) +
\chi_{G_0}u &\mbox{ in } Q,\\
\ns\ds y = 0 & \mbox{ on } \Si,\\
\ns\ds y(0) = y_0, \q y_t(0) = y_1 & \mbox { in } G.
\end{array}
\right.
\end{eqnarray}
For system (\ref{system3}), the state variable and control variable,
state space and control space, controllability, are chosen/defined
in a similar way as that for system (\ref{we}). Concerning the
nonlinearity $h(\cd)$, we need the following assumption
(\cite{DZZ}).

\begin{assum} \label{condition of d5}
Function $h(\cdot) \in C(\mathbb{R})$ is locally
Lipschitz-continuous, and for some $r\in [0,{\frac{3}{2}})$, it
satisfies that
\begin{eqnarray}\label{h}
\lim_{|s|\to\infty}\frac{\int_0^1 h_s (\tau s)d\tau}{\ln^r(1+ |s|)}
= 0.
\end{eqnarray}
\end{assum}

As mentioned in \cite{DZZ}, proceeding as in the proof of
\cite[Theorem 2.2]{Fu-Yong-Zhang}, i.e., by the linearization
approach (thanks to the second conclusion in Theorem
\ref{oTheorem5}), noting that the embedding $H_0^1(G)\hookrightarrow
L^2(G)$ is compact, and using the fixed point technique, one can
show the following result.

\begin{theorem}
Assume that Conditions \ref{condition of d-1}, \ref{condition of d}
and \ref{condition of d5} are satisfied, and $T$ and $G_0$ are given
as in Theorem \ref{oTheorem5}. Then system (\ref{slheat1}) is
exactly controllable.
\end{theorem}

Due to the blow-up and the finite propagation speed of solutions to
hyperbolic equations, one cannot expect exact controllability of
system (\ref{slheat1}) for nonlinearities of the form (\ref{h}) with
$r>2$. One could expect the system to be controllable for $r\leq2$.
However, in view of Remark \ref{ICMr1}, the usual fixed point method
cannot be applied for $r\ge3/2$. Therefore, when $n\ge 2$, the
controllability problem for system (\ref{system3}) is open for
$3/2\le r\leq2$.

\begin{remark}
Note that the above ``$3/2$ logarithmic growth" phenomenon (arising
in the global exact controllability for nonlinear PDEs) does not
occur in the pure PDE problem, and therefore the study of nonlinear
controllability is of independent interest. More precisely, this
means that for the controllability problem of nonlinear systems,
there exist some extra difficulties.
\end{remark}

\subsection{Controllability of Quasilinear PDEs}

In this subsection,  we consider the controllability of quasilinear
parabolic/hyperbolic equations.

We begin with the following controlled quasilinear hyperbolic
equation:
\begin{eqnarray}\label{qhsystem1}
\left\{
\begin{array}{lll}
\ds y_{tt} - \sum_{i,j=1}^n (h^{ij}y_{x_i})_{x_j} =
F(t,x,y,\nabla_{t,x}y,\nabla^2_{t,x}y) +qy+ \phi_{G_0}u, &\mbox{ in
}
Q,\\
\ns\ds y = 0, &\mbox{ on } \Si,\\
\ns\ds y(0) = y_0, \q y_t(0) = y_1, &\mbox{ in } G.
\end{array}
\right.
\end{eqnarray}
Here, $\big(h^{i j}\big)_{1\le i,j\le n}\in H^{s+1}(G;
\mathbb{R}^{n\times n})$ and  $\ q\in H^{s}(Q)$ with $\ds s
> \frac{n}{2}+1$, and similar to \cite{ZL}, the nonlinear term
$F(\cdot)$ has the form
$$
F(t,x,y,\nabla_{t,x}y,\nabla^2_{t,x}y) = \sum_{i=1}^n \sum_{\a=0}^n
f_{i\a}(t,x,\nabla_{t,x}y)y_{x_i x_\a} + O(|y|^2 +
|\nabla_{t,x}y|^2),
$$
where $f_{i\a}(t,x,0)=0$ and $x_0 = t$, $\phi_{G_0}$ is a
nonnegative smooth function defined on $\overline{G}$ and satisfying
$\ds\min_{x\in\cl{G_0}}\phi(x)>0$. In system (\ref{qhsystem1}), as
before, $(y,y_t)$ is the {\it state variable} and $u$ is the {\it
control variable}. However, as we shall see later, the state space
and control space have to be chosen in a different way from those
used in the linear/semilinear setting.

The controllability of quasilinear hyperbolic equations is well
understood in one space dimension (\cite{TTLi}). With regard to the
multidimensional case, we introduce the following assumption.

\begin{assum}
\label{condition of d6}
The linear part in (\ref{qhsystem1}), i.e., hyperbolic equation
\begin{eqnarray}
\left\{
\begin{array}{lll}
\ds y_{tt} - \sum_{i,j=1}^n (h^{ij}y_{x_i})_{x_j} =qy+
\chi_{G_0}u, &\mbox{ in } Q,\\
\ns\ds y=0, &\mbox{ in }\Si,\\
\ns\ds y(0)=y_0,\ y_t(0)=y_1, &\mbox{ in }G
\end{array}
\right.
\end{eqnarray}
is exactly controllable in $H_0^1(G)\t L^2(G)$ at some time $T$.
\end{assum}

Theorem \ref{oTheorem5} provides a sufficient condition to guarantee
Condition \ref{condition of d6} is satisfied. The following result
is a slight generalization of that shown in \cite{35}.

\begin{theorem}\label{qth1}
Assume Condition \ref{condition of d6} holds. Then, there is a
sufficiently small $\e_0>0$ such that for any $(y_0,y_1),
(z_0,z_1)\in \big(H^{s+1}(G)\bigcap H^1_0(G)\big)\t H^s(G)$
satisfying $|(y_0,y_1)|_{ H^{s+1}(G)\t H^s(G)}<\e_0$, $|(z_0,z_1)|_{
H^{s+1}(G)\t H^s(G)}<\e_0$ and the compatibility condition, one can
find a control $\ds u\in \bigcap_{k=0}^{s-2}C^k([0,T];H^{s-k}(G)$
such that the corresponding solution of system (\ref{qhsystem1})
verifies $y(T)=z_0$ and $y_t(T)=z_1$ in $G$.
\end{theorem}

The key in the proof of Theorem \ref{qth1} is to reduce the local
exact controllability of quasilinear equations to the exact
controllability of the linear equation by means of a new unbounded
perturbation technique (developed in \cite{35}), which is a
universal approach to solve the local controllability problem for a
large class of quasilinear time-reversible evolution equations.

Note however that the above approach does not apply to the
controllability problem for quasilinear time-irreversible evolution
equations, such as the following controlled quasilinear parabolic
equation:
\begin{eqnarray}\label{system1}
\left\{
\begin{array}{lll}
y_t-\displaystyle\sum^{n}_{i,j=1}(a^{i j}(y) y_{x_i})_{x_j}=\chi_{G_0} u  \ \ \ \ &\mbox{ in }Q,\\
 \ns
y=0 &\mbox{ on }\Sigma,\\
 \ns
y(0)=y_0 &\mbox{ in }G.
\end{array}
\right.
\end{eqnarray}
In (\ref{system1}),  $y$ is the {\it state variable} and $u$ is the
{\it control variable}, the nonlinear matrix-valued function
$\big(a^{i j}\big)_{1\le i,j\le n}\in C^2(\mathbb{R};
\mathbb{R}^{n\times n})$  is locally positive definite. One can find
very limited papers on the controllability of quasilinear
parabolic-type equations (\cite{23} and the references therein). One
of the main difficulty to solve this problem is to show the ``good
enough" regularity for solutions of system (\ref{system1}) with a
desired control.

We introduce the dual system of the linearized equation of
(\ref{system1}).
\begin{eqnarray}\label{lsystem1}
 \left\{
\begin{array}{lll}
p_t-\displaystyle\sum^{n}_{i,j=1}(p^{i j} p_{x_i})_{x_j}=0  \ \ \ &\mbox{ in }Q,\\
 \ns
p=0 &\mbox{ on }\Sigma,\\
 \ns
 p(0)=p_0 &\mbox{ in }G,
\end{array}
\right.
\end{eqnarray}
where $\big(p^{ij}\big)_{1\le i,j\le n}$ is assumed to satisfy
Condition \ref{condition of d-0}. Put $ B=1+\sum\limits^{n}_{i,j=1}
|p^{ij}|^2_{C^1(\overline{Q})}$. Starting from Theorem
\ref{uTheorem}, one can show the following observability result
(\cite{23}).

\begin{theorem}\label{oTheorem1}
There exist suitable real functions $\alpha$ and $\varphi$, and a
constant $C_0=C_0(\rho, n, G, T)>0$, such that for any $\lambda\geq
C_0
 e^{C_0B}$, solutions of (\ref{lsystem1}) satisfy
 \bel{opw}
| p(T)|_{L^2(G)} \leq C e^{e^{C B}}
\big|e^{\lambda\alpha}|\varphi|^{3/2}
p\big|_{L^2(Q_{G_0})},\qq\forall\; p_0\in L^2(G).
 \ee
\end{theorem}

In Theorem \ref{oTheorem1}, the observability constant in
(\ref{opw}) is obtained explicitly in the form of $C e^{e^{C B}}$ in
terms of the $C^1$-norms of the coefficients in the principal
operator appeared in the first equation of (\ref{lsystem1}). This is
the key in the argument of fixed point technique to show the
following local controllability of system (\ref{system1})
(\cite{23}).

\begin{theorem}\label{Theorem1}
There is a constant $\gamma>0$ such that, for any initial value
$y_0\in C^{2+\frac{1}{2}}(\overline{G})$ satisfying
$|y_0|_{C^{2+\frac{1}{2}}(\overline{G})}\leq \gamma$ and the first
order compatibility condition, one can find a control $u\in
C^{\frac{1}{2}, \frac{1}{4}}(\overline{Q})$ with $\supp
u\subseteq[0,T]\times G_0$ such that the solution $y$ of system
(\ref{system1}) satisfies  $y(T)=0$ in $G$.
\end{theorem}

From Theorem \ref{Theorem1}, it is easy to see that the {\it state
space} and {\it control space} for system (\ref{system1}) are chosen
to be $C^{2+\frac{1}{2}}(\overline{G})$ and  $C^{\frac{1}{2},
\frac{1}{4}}(\overline{Q})$, respectively. The key observation in
\cite{23} is that, thanks to an idea in \cite{barbu}, for smooth
initial data, the regularity of the null-control function for the
linearized system can be improved, and therefore, the fixed point
method is applicable.

\subsection{Stabilization of hyperbolic equations and further comments}

In this subsection, we give more applications of Theorem
\ref{uTheorem} to the stabilization of hyperbolic equations and
comment other applications of this theorem and some related open
problems.

One of the main motivation to introduce the
controllability/obseervability theory is to design the feedback
regulator (\cite{15}). Stimulated by \cite{Lions}, there exist a lot
of works addressing the stabilization problem of PDEs from the
control point of view. To begin with, we fix a nonnegative function
$a\in L^{\i}(\G)$ such that
 $\big\{x\in\G\;\big|\;a(x)>0\big\}\neq\emptyset$,
and consider the following hyperbolic equation with a boundary
damping:
 \bel{0a1}\left\{\ba{ll}\ds
 u_{tt}-\sum_{j,k=1}^n(h^{jk}u_{x_j})_{x_k}=0 &\hb{ in } (0,\infty)\t G,\\
 \ns\ds
\sum_{j,k=1}^nh^{jk}u_{x_j}\nu_k+a(x)u_t=0&\hb{ on } (0,\infty)\t \G,\\
\ns\ds u(0)=u^0, \q u_t(0)=u^1&\hb{ in } G.
 \ea\right.\ee

Put $
 H\=\left\{(f,g)\in H^1(G)\t L^2(G)\;\left|\;\int_G fdx=0\right\}\right.$, which is a Hilbert
space with the canonic norm. Define an unbounded operator $\cA: \;
H\to H$ by
 $$\left\{\ba{ll}
 \3n \cA\=\left(\ba{cc}
 0&I\\ \ns\ds\sum_{j,k=1}^n\pa_{x_k}(h^{jk}\pa_{x_j})&0
 \ea
 \right),\\
 \ns\ds \3n D(\cA)\=\Big\{u=(u^0,u^1)\in H\;\Big|\  \cA u\in H,\  \(\sum_{j,k=1}^nh^{jk}u^0_{x_j}\nu_k+au^1\)\|_{\G}=0\Big\}. \ea\right.
 $$
It is easy to show that $\cA$ generates an $C_0$-semigroup
$\{e^{t\cA}\}_{t\in \mathbb{R}}$ on $H$. Hence, system (\ref{0a1})
is well-posed in $H$. Clearly, $H$ is the {\it finite energy space}
of system (\ref{0a1}). One can show that the energy of any solution
of (\ref{0a1}) tends to zero as $t\to\infty$ (There is no any
geometric conditions on $\G$).

Starting from Theorem \ref{uTheorem}, one can show the following
result, which is a slight improvement of the main result in
\cite{Fu2}:
 \bt\label{0t1}
Assume Conditions \ref{condition of d-1} holds. Then solutions $u\in
C([0,\infty);$ $ D(\cA))\bigcap C^1([0,\infty);\;H)$ of system
(\ref{0a1}) satisfy
  \bel{0a3}
 ||(u,u_t)||_{H}\le \frac{C}{ \ln (2+t)}||(u^0,u^1)||_{D(\cA)},\q\forall\;(u^0,u^1)\in D(\cA),\ \forall\; t>0.
  \ee
 \et

Next, we consider a semilinear hyperbolic equation with a local
damping:
 \bel{icm0a1}\left\{\ba{ll}\ds
 u_{tt}-\sum_{j,k=1}^n(h^{jk}u_{x_j})_{x_k}+f(u)+b(x)g(u_t,\nabla u)=0 &\hb{ in } (0,\infty)\t G,\\
 \ns\ds
u=0&\hb{ on } (0,\infty)\t \G,\\
\ns\ds u(0)=u^0, \q u_t(0)=u^1&\hb{ in } G.
 \ea\right.\ee
In (\ref{icm0a1}), $h^{jk}$ is supposed to satisfy Conditions
\ref{condition of d-1} and \ref{condition of d}; $f:\mathbb{R}\to
\mathbb{R}$ is a differentiable function satisfying $f(0) = 0$,
$sf(s)\ge  0$ and $|f'(s)| \le C (1 + |s|^q)$ for any
$s\in\mathbb{R}$, where $q \ge 0$ and $(n -2)q \le 2$; $b$ is a
nonnegative function satisfying $\ds\min_{x\in G_0} b(x)>0$, where
$G_0$ is given in Theorem \ref{oTheorem5}; and $g :
\mathbb{R}^{n+1}\to \mathbb{R}$ is a globally Lipschitz function
satisfying $g(0, w) = 0$,  $|g(r,w) - g(r_1, w_1|\le C (|r - r_1| +
|w - w_1|)$ and $g(r,w)r \ge c_0r^2$ for some $c_0
> 0 $, any $w, w_1\in\mathbb{R}^n$ and any $r, r_1\in
\mathbb{R}$.

Define the energy of any solution $u$ to (\ref{icm0a1}) by setting
 $$
 E(t)=\frac12\int_G\Big[|u_t|^2+\sum_{j,k=1}^nh^{jk}u_{x_j}u_{x_k}\Big]dx+\int_G\int_0^uf(s)dsdx.
 $$
Starting from Theorem \ref{uTheorem}, one can show the following
stabilization result for system (\ref{icm0a1}) (\cite{Tebou}).

\bt
Let $(u^0,u^1)\in H_0^1(G)\times L^2(G)$. Then there exist positive
constants $M$ and $r$, possibly depending on $E(0)$, such that the
energy $E(t)$ of the solution of (\ref{icm0a1}) satisfies $E(t)\le
Me^{-rt}E(0)$ for any $t\ge 0$.
\et

Several comments are in order.

\begin{remark}\label{icmrem2}
In \cite{LR}, the authors need $C^\infty$-regularity for the data to
establish Theorem \ref{icm2t1}. Recently, based on Theorem
\ref{uTheorem}, this result was extended in \cite{25} as follows:
Denote by $\{\l_i\}^\infty_{i=1}$ the eigenvalues of any general
elliptic operator of second order (with $C^{1}$-principal part
coefficients) on $\O$ (of class $C^{2}$) with Dirichlet or Robin
boundary condition, and $\{e_i\}^{\infty}_{i=1}$ the corresponding
eigenvectors satisfying $|e_i|_{L^2(\O)} = 1$. Then, for any $r
> 0$, it holds
 $$
\sum_{\l_{i} \leq r}|a_{i}|^{2} \leq C
e^{C\sqrt{r}}{\int_{G_0}}\bigg| \sum_{\l_{i} \leq r}a_{i}e_{i}(x)
\bigg|^2dx,\ \ \forall\, \{a_{i}\}_{\l_{i} \leq r }\hb{ with }a_{i}
\in \mathbb{C}.
 $$
\end{remark}

\begin{remark}
As indicated in \cite{Fu1,LTZ,LTZ2}, Theorem \ref{uTheorem} can be
employed to study the global unique continuation and inverse
problems for some PDEs. Note also that this Carleman estimate based
approach can be applied to solve some optimal control problems
(\cite{WW}).
\end{remark}

\begin{remark}
In practice, constrained controllability is more realizable. It is
shown in \cite{[40]} that the study of this problem is unexpectedly
difficult even for the $1-d$ wave equation and heat equation. We
refer to \cite{22} for an interesting example showing that this
problem is nontrivial even if the control is effective everywhere in
the domain in which the system is evolved.
\end{remark}

\begin{remark}
Note that the above mentioned approach applies mainly to the
controllability, observability and stabilization of second order
non-degenerate PDEs. It is quite interesting to extend it to the
coupled and/or higher order systems, or degenerate systems but in
general, this is nontrivial even for linear problems
(\cite{CMV,ZZ1}).
\end{remark}

\begin{remark}
Similar to other nonlinear problems, nonlinear controllability
problems are usually quite difficult. It seems that there is no
satisfactory controllability results published for nonlinear
hyperbolic-parabolic coupled equations. Also, there exists no
controllability results for fully nonlinear PDEs. In the general
case, of course, one could expect only local results. Therefore, the
following three problems deserve deep studies: 1) The
characterization of the controllability subspace; 2) Controllability
problem with (sharp) lower regularity for the data; 3) The problem
that cannot be linearized. Of course, all of these problems are
usually challenging.
\end{remark}

\section{The stochastic case}\label{s4}

In this section, we extend some of the results/approaches in Section
\ref{s3} to the stochastic case. As we shall see later, the
stochastic counterpart is far from satisfactory, compared to the
deterministic setting.

In what follows, we fix a complete filtered probability space
$(\O,\cF,\{\cF_t\}_{t\ge 0},P)$ on which a one dimensional standard
Brownian motion $\{B(t)\}_{t\ge0}$ is defined. Let $H$ be a
Fr\'echet space. Denote by $L_{\cF}^2(0,T;H)$ the Fr\'echet space
consisting of all $H$-valued $\{\cF_t\}_{t\ge 0}$-adapted processes
$X(\cd)$ such that $\mathbb{E}(|X(\cd)|_{L^2(0,T;H)}^2)<\i$,  with
the canonical quasi-norms; by $L_{\cF}^\i(0,T;H)$ the Fr\'echet
space consisting of all $H$-valued $\{\cF_t\}_{t\ge 0}$-adapted
bounded processes, with the canonical quasi-norms; and by
$L_{\cF}^2(\O; C([0,T];H))$ the Fr\'echet space consisting of all
$H$-valued $\{\cF_t\}_{t\ge 0}$-adapted continuous processes
$X(\cd)$ such that $ \mathbb{E}(|X(\cd)|_{C([0,T];H)}^2)<\i$, with
the canonical quasi-norms (similarly, one can define $L^{2}_{\cal
F}(\O;C^{k}([0,T];H))$ for $k\in \mathbb{N}$).

\subsection{Stochastic Parabolic  Equations}

We begin with the following stochastic parabolic equation:
 \bel{hh6.1}
 \left\{
 \ba{ll}
 \ds dz-\sum_{i,j=1}^n(p^{ij}z_{x_i})_{x_j}dt=[\lan a,\n z\ran+bz]dt+cz
 dB(t)&\hb{ in }Q,\\
 \ns
 z=0&\hb{ on }\Si,\\
 \ns
 z(0)=z_0&\hb{ in }G
 \ea
 \right.
 \ee
with suitable coefficients $a,b$ and $c$, where $p^{ij}\in
C^2(\cl{Q})$ is assumed to satisfy Condition \ref{condition of d-0}
(Note that, technically we need here more regularity for $p^{ij}$
than the deterministic case). We are concerned with an observability
estimate for system (\ref{hh6.1}), i.e., to find a constant
$\cC=\cC(a,b,c,T)>0$ such that solutions of (\ref{hh6.1}) satisfy
 \bel{12e2}
 |z(T)|_{L^2(\O,\cF_T,P;L^2(G))}\le \cC|z|_{L^2_{\cF}(0,T;L^2(G_0))},\q\forall\;z_0\in L^2(\O,\cF_0,P;L^2(G)).
 \ee

Similar to Theorem \ref{uTheorem}, we have the following weighted
identity (\cite{Tang-Zhang1}).

\bt\label{icmc1t1}
Let $m\in\mathbb{N}$, $b^{ij}=b^{ji}\in
L_{\cF}^2(\O;C^1([0,T];W^{2,\infty}(\mathbb{R}^m)))$
($i,j=1,2,\cdots,m$), $\ell\in C^{1,3}((0,T)\t\mathbb{R}^m)$ and
$\Psi\in C^{1,2}((0,T)\t\mathbb{R}^m)$. Assume $u$ is an
$H^2(\mathbb{R}^m)$-valued continuous semi-martingale. Set
 $\th=e^{\ell }$ and $v=\th u$. Then for $\ae\2n$ $x\in \mathbb{R}^m$ and $P$-$\as\2n$
 $\o\in\O$,
 $$\ba{ll}
 \displaystyle
 2\int_0^T\th\[-\sum_{i,j=1}^m (b^{ij}v_{x_i})_{x_j}+Av\]\[du-\sum_{i,j=1}^m(b^{ij}u_{x_i})_{x_j}dt\]+2\int_0^T\sum_{i,j=1}^m(b^{ij}v_{x_i}dv)_{x_j}\\
 \noalign{\ss}
 \displaystyle\q+ 2\int_0^T\sum_{i,j=1}^m\[\sum_{i',j'=1}^m\(2b^{ij} b^{i'j'}\ell_{x_{i'}}v_{x_i}v_{x_{j'}}
 -b^{ij}b^{i'j'}\ell_{x_i}v_{x_{i'}}v_{x_{j'}}\)\\
 \noalign{\ss}
 \displaystyle\q+\Psi
 b^{ij}v_{x_i}v- b^{ij}\(A\ell_{x_i}+\frac{\Psi_{x_i}}{2}\)v^2\]_{x_j}dt\\
 \noalign{\ss}
 \displaystyle
 =2\int_0^T\sum_{i,j=1}^m \Big\{\sum_{i',j'=1}^m\[2b^{ij'}\(b^{i'j}\ell_{x_{i'}}\)_{x_{j'}} -
 \(b^{ij}b^{i'j'}\ell_{x_{i'}}\)_{x_{j'}}\]-\frac{b_t^{ij}}{2}+\Psi b^{ij}
 \Big\}v_{x_i}v_{x_j}dt\\
 \noalign{\ss}
 \displaystyle\q
 +\int_0^TBv^2dt+2\int_0^T\[-\sum_{i,j=1}^m (b^{ij}v_{x_i})_{x_j}+Av\]\[-\sum_{i,j=1}^m (b^{ij}v_{x_i})_{x_j}+(A-\ell_t)v\]dt\\
 \noalign{\ss}
 \displaystyle\q+\(\sum_{i,j=1}^m
 b^{ij}v_{x_i}v_{x_j}+Av^2\)\Big|_0^T\\
 \noalign{\ss}
 \displaystyle\q-\int_0^T\th^2\sum_{i,j=1}^m
 b^{ij}[(du_{x_i}+\ell_{x_i}du)(du_{x_j}+\ell_{x_j}du)]-\int_0^T\th^2A(du)^2,
 \ea
 $$
where
 $$
\left\{
 \ba{ll}
 \ds A\=-\sum_{i,j=1}^m
 (b^{ij}\ell_{x_i}\ell_{x_j}-b^{ij}_{x_j}\ell_{x_i}
 -b^{ij}\ell_{x_ix_j})-\Psi,\\
  \ns
 \ds
 B\=2\[A\Psi-
 \sum_{i,j=1}^m(Ab^{ij}\ell_{x_i})_{x_j}\] -A_t-\sum_{i,j=1}^m (b^{ij}\Psi_{x_j})_{x_i}.
  \ea
\right.
$$
\et

\br
Note that, in Theorem \ref{icmc1t1}, we assume only the symmetry for
matrix $\big(b^{ij}\big)_{1\le i,j\le n}$ (without assuming the
positive definiteness). Hence, this theorem can be applied to study
not only the observability/controllability of stochastic parabolic
equations, but also similar problems for deterministic parabolic and
hyperbolic equations, as indicated in Section \ref{s3}. In this way,
we give a unified treatment of controllability/observability
problems for some stochastic and deterministic PDEs of second order.
\er

Starting from Theorem \ref{icmc1t1}, one can show the following
observability result in \cite{Tang-Zhang1} (See \cite{barbu1} and
the references therein for some earlier results).

\bt\label{1t4}
Assume that
$$a\in  L^{\i}_{\cF}(0,T;L^{\i}(G;\dbR^n)), \ \ b\in
L^{\i}_{\cF}(0,T;L^{n^*}(G)), \ \  c\in
L^{\i}_{\cF}(0,T;W^{1,\i}(G)),
 $$
where $n^*\ge 2$ if $n=1$; $n^*> 2$ if $n=2$; $n^*\ge n$ if $n\ge
3$. Then there is a constant $\cC=\cC(a,b,c,T)>0$ such that all
solutions $z$ of system (\ref{hh6.1}) satisfy (\ref{12e2}).
Moreover, the observability constant $\cC$ may be bounded as
 $$
 \cC(a,b,c,T)=Ce^{C[T^{-4}(1+\tau^2)+T\tau^2]},
 $$
with
$\tau\=|a|_{L^{\i}_{\cF}(0,T;L^{\i}(G;\dbR^n))}+|b|_{L^{\i}_{\cF}(0,T;L^{n^*}(G))}+|c|_{
L^{\i}_{\cF}(0,T;W^{1,\i}(G))}$.
\et

As a consequence of Theorem \ref{1t4}, one can deduce a
controllability result for backward stochastic parabolic equations.
Unlike the deterministic case, the study of controllability problems
for forward stochastic differential equations is much more difficult
than that for the backward ones. We refer to \cite{Peng} for some
important observation in this respect. It deserves to mention that,
as far as I know, there exists no satisfactory controllability
result published for forward stochastic parabolic equations. Note
however that, as a consequence of Theorem \ref{icm2t1} and its
generalization (see Remark \ref{icmrem2}), one can deduce a null
controllability result for forward stochastic parabolic equations
with time-invariant coefficients (\cite{25}).

Theorem \ref{icmc1t1} has another application in global unique
continuation of stochastic PDEs. To see this, we consider the
following stochastic parabolic equation:
  \begin{equation}\label{icmhh6.1}
  \cF z\equiv dz-\sum_{i,j=1}^n(f^{ij}z_{x_i})_{x_j}dt=[\lan a_1,\n z\ran+b_1z]dt+c_1z
  dB(t)\q\hb{ in }Q,
  \end{equation}
where $f^{ij}\in C^{1,2}([0,T]\times G)$ satisfy $f^{ij}=f^{ji}$
($i,j=1,2,\cdots,n$) and for any open subset $G_1$ of $G$, there is
a constant $s_0=s_0(G_1)>0$ so that
  $\ds\sum_{i,j=1}^nf^{ij}\xi^{i}\xi^{j}
  \geq s_0|\xi|^{2}$
for all $(t,x,\xi)\equiv (t,x,\xi^{1},\cdots, \xi^{n})
\in(0,T)\times G_1\t \dbR^{n}$; $a_1\in
L^{\i}_{\cF}(0,T;L_{loc}^{\i}(G;\dbR^n))$, $b_1\in
L^{\i}_{\cF}(0,T;L_{loc}^{\i}(G))$, and $c_1\in L^{\i}_{\cF}(0,T;$
$W_{loc}^{1,\i}(G))$.

Starting from Theorem \ref{icmc1t1}, one can show the following
result (\cite{[39]}).

\bt\label{1t3} Any solution $ z\in
L_{\cF}^2(\O;C([0,T];L_{loc}^2(G))) \bigcap $ $
L_{\cF}^2(0,T;H_{loc}^1(G)) $ of $(\ref{icmhh6.1})$ vanishes
identically in $Q\t\O, \as dP$ provided that $z=0$ in $Q_{G_0}\t\O,
\as dP$.
\et

Note that the solution of a stochastic equation is generally {\it
non-analytic in time} even if all coefficients of the equation are
constants. Therefore, one cannot expect a Holmgren-type uniqueness
theorem for stochastic equations except for some very special cases.
On the other hand, the usual approach to employ Carleman-type
estimate for the unique continuation needs to localize the problem.
The difficulty of our present stochastic problem consists in the
fact that {\it one cannot simply localize the problem as usual}
because the usual localization technique may change the adaptedness
of solutions, which is a key feature in the stochastic setting. In
equation (\ref{icmhh6.1}), for the space variable $x$, we may
proceed as in the classical argument. However, for the time variable
$t$, due to the adaptedness requirement, we will have to treat it
separately and globally. We need to introduce {\it partial global
Carleman estimate} (indeed, global in time) even for local unique
continuation for stochastic parabolic equation. Note that this idea
comes from the study of controllability problem even though unique
continuation itself is purely an PDE problem.

\subsection{Stochastic Hyperbolic Equations}

We consider now the following stochastic wave equation:
\begin{eqnarray}\label{ssystem3}
\left\{
\begin{array}{lll}\ds
\ds dz_{t} - \D z dt =(a_1 z_t + \lan a_2,\nabla z\ran + a_3 z + f
)dt + (a_4 z + g)dB(t) & {\mbox { in }} Q,
 \\
\ns\ds  z = 0 & \mbox{ on } \Si, \\
\ns\ds  z(0) = z_0, \q z_{t}(0) = z_1 & \mbox{ in } G,
\end{array}
\right.
\end{eqnarray}
where $a_1\in L^{\infty}_{\cF}(0,T;L^{\infty}(G))$, $a_2\in
 L^{\infty}_{\cF}(0,T;L^{\infty}(G;\mathbb{R}^n))$, $a_3\in L^{\infty}_{\cF}(0,T;L^n(G))$,
 $a_4\in
 L^{\infty}_{\cF}(0,T;L^\infty(G))$, $f\in L^2_{\cF}(0,T;L^2(G))$, $g\in L^2_{\cF}(0,T;L^2(G))$ and
 $(z_0, z_1) \in L^2(\O,{\cal F}_0, P; H^1_0(G) \t L^2(G))$. We shall derive an observability estimate for
(\ref{icmhh6.1}), i.e., find a constant $\cC(a_1,a_2,a_3,a_4)>0$
such that solutions of system (\ref{icmhh6.1}) satisfy
 \bel{icm12e2}
 \ba{ll}\ds
 |(y(T),y_t(T))|_{L^2(\O,\cF_T,P;H_0^1(G)\times L^2(G))}\\ \ns
 \ds\le \1n\cC(a_1,a_2,a_3,a_4)\2n\left[\left|\frac{\pa y}{\pa\nu}\right|_{ L^2_{\cF}(0,T;L^2(\G_0))}\2n+|f|_{L^2_{\cF}(0,T;L^2(G))}+|g|_{
L^2_{\cF}(0,T;L^2(G))}\2n\right]\1n,\\ \ns
 \ds\qq\qq\qq\qq\forall\;(y_0,y_1)\in L^2(\O,\cF_0,P;H_0^1(G)\times
 L^2(G)).
 \ea
 \ee
where $\G_0$ is given by (\ref{boundary}) for some
$x_0\in\mathbb{R}^d\setminus \cl G$.

It is clear that, $\ds 0<R_0\=\min_{x\in G}|x-x_0|< R_1\=\max_{x\in
G}|x-x_0|$. We choose a sufficiently small constant $c\in (0,1)$ so
that $\frac{(4+5c)R_0^2}{9c}>R_1^2$. In what follows, we take $T$
sufficiently large such that
$\frac{4(4+5c)R_0^2}{9c}>c^2T^2>4R_1^2$. Our observability estimate
for system (\ref{icmhh6.1}) is stated as follows (\cite{Zhangxu3}).

\bt \label{t1} Solutions of system $(\ref{ssystem3})$ satisfy
$(\ref{icm12e2})$ with
 $$
 \ba{ll}
 \cC(a_1,a_2,a_3,a_4)\\
 \ns
 \ds=\2n C\exp\1n\left\{\1n C\1n \left[
|(a_1,a_4)|_{L^{\infty}_{\cF}(0,T;(L^\infty(G))^2)}^2+| a_2|_{
 L^{\infty}_{\cF}(0,T;L^{\infty}(G;\mathbb{R}^n))}^2  + |a_3|_{L^{\infty}_{\cF}(0,T;L^n(G))}^2\1n\right]\1n\right\}.
 \ea
$$
 \et

Surprisingly, Theorem \ref{t1} was improved in \cite{25} by
replacing the left hand side of (\ref{icm12e2}) by $
|(y_0,y_1)|_{L^2(\O,\cF_0,P;H_0^1(G)\times L^2(G))}$, exactly in a
way of the deterministic setting. This is highly nontrivial by
considering the very fact that the stochastic wave equation is
time-irreversible.

The proof of Theorem \ref{t1} (and its improvement in \cite{25}) is
based on the following identity for a stochastic hyperbolic-like
operator, which is in the spirit of Theorems \ref{uTheorem} and
\ref{icmc1t1}.

\bt \label{c1t1} Let $b^{ij}\in C^1((0,T)\t\mathbb{R}^n)$ satisfy
 $b^{ij}=b^{ji}$ ($i,j=1,2,\ldots,n$), $\ell,\ \Psi\in
C^2((0,T)\t\mathbb{R}^n)$. Assume $u$ is an
$H^2_{loc}(\mathbb{R}^n)$-valued $\{\cF_t\}_{t\ge 0}$-adapted
process such that $u_t$ is an $L^2_{loc}(\mathbb{R}^n)$-valued
semimartingale. Set
 $\th=e^{\ell }$ and $v=\th u$. Then, for a.e. $x\in \mathbb{R}^n$ and
 $P$-a.s. $\o\in\O$,
  $$
  \ba{ll}
 \displaystyle
 \th\Big(-2\ell_tv_t+2\sum_{i,j=1}^nb^{ij}\ell_{x_i}v_{x_j}+\Psi v\Big)\Big[du_t-\sum_{i,j=1}^n(b^{ij}u_{x_i})_{x_j}dt\Big]\\
\ns
 \displaystyle\q+
 \sum_{i,j=1}^n\biggl[\sum_{i',j'=1}^n\left(2b^{ij} b^{i'j'}\ell_{x_{i'}}v_{x_i}v_{x_{j'}}
 -b^{ij}b^{i'j'}\ell_{x_i}v_{x_{i'}}v_{x_{j'}}\right)-2b^{ij}\ell_tv_{x_i}v_t
 +b^{ij}\ell_{x_i}v_t^2\\
\ns
 \displaystyle\qq\qq\qq+\,\Psi b^{ij}v_{x_i}v-\Big(A\ell_{x_i}+\frac{\Psi_{x_i}}{2}\Big)b^{ij} v^2\biggr]_{x_j}dt\\
\ns
 \displaystyle\q+\,d\Big[\sum_{i,j=1}^n
 b^{ij}\ell_tv_{x_i}v_{x_j}-2\sum_{i,j=1}^nb^{ij}\ell_{x_i}v_{x_j}v_t+\ell_tv_t^2-\Psi
 v_tv+\left(A\ell_t+\frac{\Psi_t}{2}\right)v^2\Big]\\
\ns
 \displaystyle
 =\Big\{\Big[\ell_{tt}+\sum_{i,j=1}^n(b^{ij}
 \ell_{x_i})_{x_j}-\Psi\Big]v_t^2-2\sum_{i,j=1}^n\left[(b^{ij}\ell_{x_j})_t+b^{ij}\ell_{tx_j}\right]v_{x_i}v_t\\
\ns
 \displaystyle\q+\sum_{i,j=1}^n\Big[(b^{ij}\ell_t)_t+\sum_{i',j'=1}^n\left(2b^{ij'}(b^{i'j}\ell_{x_{i'}})_{x_{j'}} -
 (b^{ij}b^{i'j'}\ell_{x_{i'}})_{x_{j'}}\right)+\Psi b^{ij}
 \Big]v_{x_i}v_{x_j}\\
\ns
 \displaystyle\q
 +\,{\it Bv}^2+\Big(-2\ell_tv_t+2\sum_{i,j=1}^nb^{ij}\ell_{x_i}v_{x_j}+\Psi v\Big)^2\Big\}dt+\th^2\ell_t(du_t)^2,
 \ea
 $$
where $(du_t)^2$ denotes the quadratic variation process of $u_t$,
 $$
\left\{
 \ba{ll}
 \ds A\=(\ell_t^2-\ell_{tt})-\sum_{i,j=1}^n
 (b^{ij}\ell_{x_i}\ell_{x_j}-b^{ij}_{x_j}\ell_{x_i}
 -b^{ij}\ell_{x_ix_j})-\Psi,\\
\ns
 \ds
 B\=A\Psi+(A\ell_t)_t-
 \sum_{i,j=1}^n(Ab^{ij}\ell_{x_i})_{x_j} +\frac{1}{2}\Big[\Psi_{tt}-\sum_{i,j=1}^n (b^{ij}\Psi_{x_i})_{x_j}\Big].
 \ea
\right.
 $$
\et

\subsection{Further comments}

Compared to the deterministic case, the
controllability/observability of stochastic PDEs is in its ``enfant"
stage. Therefore, the main concern of the
controllability/observability theory in the near future should be
that for stochastic PDEs. Some most relevant open problems are
listed below.

\begin{enumerate}

\item[$\bullet$] {\bf Controllability of forward stochastic PDEs}.
Very little is known although there are some significant progress in
the recent work \cite{25}. Also, it would be quite interesting to
extend the result in \cite{BLR} to the stochastic setting but this
seems to be highly nontrivial.

\item[$\bullet$]{\bf Controllability of nonlinear stochastic PDEs}. Almost nothing is known in this direction although
there are some papers addressing the problem in abstract setting by
imposing some assumption which is usually very difficult to check
for the nontrivial case.

\item[$\bullet$]{\bf Stabilization and inverse problems for stochastic PDEs}. Almost nothing is known in this respect.

\end{enumerate}

\section*{Acknowledgement}

This paper is a summary of part of the work I have done in close
collaboration with my colleagues, coworkers and students. I am
grateful to all of them. In particular, I would like to express my
sincerely gratitude to my mentors, Xunjing~Li, Jiongmin~Yong and
Enrique~Zuazua, for their continuous encouragement and for so many
fruitful discussions and suggestions that have been extremely
influential on the formulation and solution of most of problems
addressed here. Also, I would like to mention some of my colleagues
with who I had the opportunity to develop part of the theory and
learn so many things and, in particular, Xiaoyu~Fu, Xu~Liu, Qi~L\"u,
Kim~Dang~Phung, Shanjian~Tang, Gengsheng~Wang, Jiongmin~Yong and
Enrique~Zuazua. Finally, I thank Viorel~Barbu, Piermarco~Cannarsa,
Xiaoyu~Fu, Yamilet~Quintana, Louis~T\'ebou, Marius~Tucsnak and
Gengsheng~Wang for their useful comments on the first version of
this paper that allowed to improve its presentation and to avoid
some inaccuracies.


\end{document}